\documentclass[12pt,reqno]{amsart}
\usepackage{enumerate, latexsym, amsmath, amsfonts, amssymb, amsthm, color, hyperref}
\def\pmod #1{\ ({\rm{mod}}\ #1)}
\def\Z{\Bbb Z}
\def\N{\Bbb N}

\def\l{\left}
\def\r{\right}
\def\bg{\bigg}
\def\({\bg(}
\def\){\bg)}
\def\t{\text}
\def\f{\frac}

\def\mo{{\rm{mod}\ }}

\def\sm{\setminus}
\def\bi{\binom}

\def\eq{\equiv}

\def\Proof{\noindent{\it Proof}}
\def\Ack{\medskip\noindent {\bf Acknowledgments}}
\theoremstyle{plain}
\newtheorem{theorem}{Theorem}

\newtheorem{lemma}{Lemma}

\theoremstyle{definition}

\theoremstyle{remark}
\newtheorem{remark}{Remark}

 \vspace{4mm}

\begin{document}

\hbox{Accepted by Linear and Multilinear Algebra.}
\medskip

\title
[{Proof of three conjectures on determinants}]
{Proof of three conjectures on determinants related to \\quadratic residues}

\author
[Darij Grinberg, Zhi-Wei Sun and Lilu Zhao] {Darij Grinberg, Zhi-Wei Sun and Lilu Zhao}

\address {(Darij Grinberg) Mathematics Department, Drexel University,
 Philadelphia, PA, USA}
\email{darijgrinberg@gmail.com}

\address{(Zhi-Wei Sun) Department of Mathematics, Nanjing
University, Nanjing 210093, People's Republic of China}
\email{zwsun@nju.edu.cn}

\address {(Lilu Zhao) School of Mathematics, Shandong University, Jinan 250100, People's Republic of China}
\email{zhaolilu@sdu.edu.cn}

\keywords{Determinant, divisibility, Jacobi symbol, Vandermonde-type determinant.
\newline \indent The second author is the corresponding author, and supported by the Natural Science Foundation of China (grant no. 11971222).}

\begin{abstract} In this paper we confirm three conjectures of Z.-W. Sun on determinants.
We first show that any odd integer $n>3$ divides the determinant
$$\left|(i^2+dj^2)\left(\frac{i^2+dj^2}n\right)\right|_{0\le i,j\le (n-1)/2},$$
where $d$ is any integer and $(\frac{\cdot}n)$ is the Jacobi symbol. Then we prove some divisibility results concerning $|(i+dj)^n|_{0\le i,j\le n-1}$ and $|(i^2+dj^2)^n|_{0\le i,j\le n-1}$, where $d\not=0$ and $n>2$ are integers. Finally, for any odd prime $p$ and integers $c$ and $d$ with $p\nmid cd$, we determine completely the Legendre symbol $\left(\frac{S_c(d,p)}p\right)$, where $S_c(d,p):=\left|\left(\frac{i^2+dj^2+c}p\right)\right|_{1\le i,j\le(p-1)/2}$.
\end{abstract}
\maketitle

\section{Introduction}
\setcounter{lemma}{0}
\setcounter{theorem}{0}
\setcounter{corollary}{0}
\setcounter{remark}{0}
\setcounter{equation}{0}

For an $n\times n$ matrix $[a_{ij}]_{1\le i,j\le n}$ over a commutative ring with identity, we shall denote its determinant by $|a_{ij}|_{1\le i,j\le n}$.
In this paper we study some determinants related to quadratic residues. For the standard theory of quadratic residues, one may consult \cite[Chapter 5, pp.\,50-65]{IR}.

Our first theorem in the case $d=1$ was originally conjectured by
Sun \cite[Conjecture 4.5(i)]{S19} amid a study of determinants
involving Jacobi symbols.

\begin{theorem}\label{Th1.1} Let $n>3$ be an odd integer. For any integer $d$, we have
\begin{equation}\label{0(mod n)}\bg|(i^2+dj^2)\l(\f{i^2+dj^2}n\r)\bg|_{0\le i,j\le(n-1)/2}\eq0\pmod n,
\end{equation}
where $(\f{\cdot}n)$ denotes the Jacobi symbol.
\end{theorem}
\begin{remark}\label{Rem1.1} For $n=3$ and $d\in\Z$, clearly
\begin{align*}&\bg|(i^2+dj^2)\l(\f{i^2+dj^2}n\r)\bg|_{0\le i,j\le(n-1)/2}
\\=&\left|\begin{matrix}
0&d\l(\f d3\r)\\1&(1+d)\l(\f{1+d}3\r)
\end{matrix}\right|
=-d\l(\f d3\r)\eq\begin{cases}0\pmod3&\t{if}\ 3\mid d,
\\-1\pmod 3&\t{if}\ 3\nmid d.
\end{cases}\end{align*}
\end{remark}

Let $p$ be an odd prime. R. Chapman \cite{Ch} evaluated
the determinant
$$\l|\l(\f{i+j-1}p\r)\r|_{1\le i,j\le(p+1)/2}=\l|\l(\f{i+j}p\r)\r|_{0\le i,j\le(p-1)/2},$$
and M. Vsemirnov \cite{V12,V13} determined the exact value of
$$\l|\l(\f{i-j}p\r)\r|_{0\le i,j\le(p-1)/2}$$
guessed by Chapman. Recall that $(\f ap)\eq a^{(p-1)/2}\pmod p$ for any $a\in\Z$.

Our next theorem in the case $c=0$ and $d=1$ confirms a conjecture of Sun \cite{S13} posed in 2013.

\begin{theorem}\label{Th1.2} Let $c$, $d$ and $n$ be integers with $d \neq 0$ and $n>2$.
Set
\begin{equation}\label{ab}a_n=\l|(i+dj+c)^n\r|_{0\le i,j\le n-1}\ \ \t{and}\ \ b_n=\l|(i^2+dj^2)^n\r|_{0\le i,j\le n-1}. \end{equation}
Then
\begin{equation}\label{ab'}a_n'=\f{a_nd^{-n(n-1)/2}}{(n-2)!n\prod_{k=1}^n k!} \ \t{and}\ b_n'=\f{b_nd^{-n(n-1)/2}}{2\prod_{k=1}^n(k!(2k-1)!)}
\end{equation}
are integers; in particular,
\begin{equation}\label{divide} d^{n(n-1)/2}n^2\mid a_n\ \ \t{and}\ \ d^{n(n-1)/2}(2n)!\mid b_n.
\end{equation}
Also, $(-1)^{n(n-1)/2}a_n>0$ and $(-1)^{n(n-1)/2}b_n>0$ if $d>0$ and $c\ge0$.
\end{theorem}
\begin{remark}\label{Rem1.2} Let $c,d\in\Z$ with $d\not=0$. For small values of $n$ we may compute $a_n$ via the {\tt Mathematica} command

{\tt FullSimplify[Det[Table[(i+d*j+c)${}^\wedge$n,$\{$i,0,n-1$\}$,$\{$j,0,n-1$\}$]]]}.

\noindent In particular, we get
\begin{align*}a_1=&0,\ a_2=-d (d + 2 c (1 + c + d)),
\\a_3=&-36d^3(1 + c + d) (c (2 + c) + 2 (1 + c) d),
\\a_4=&2304d^6\big(6 c (1 + c) (2 + c) (3 + c) + 18 (3 + 2 c) (1 + c (3 + c)) d
\\&+ 11 (11 + 6 c (3 + c)) d^2  +18 (3 + 2 c) d^3\big).
\end{align*}
It follows that
$$a_1'=0,\ a_2'=\f{a_2}{4d}=-\f{d+2c(1+c+d)}4,\ a_3'=\f{a_3}{36d^3}\in\Z,\ a_4'=\f{a_4}{2304d^6}\in\Z.$$
Similarly, we may compute $b_n$ via the
{\tt Mathematica} command

{\tt Expand[Det[Table[(i${}^\wedge$2+d*j${}^\wedge$2)${}^\wedge$n,$\{$i,0,n-1$\}$,$\{$j,0,n-1$\}$]]]}.

\noindent To get $b_n'$ from $b[n]=b_n$, we use the {\tt Mathematica} command

{\tt Expand[b[n]d${}\wedge$(-n(n-1)/2)/(2*Product[k!(2k-1)!,{k,1,n}])]}.

\noindent In particular, we obtain
\begin{align*}b_1'=&0,\ b_2'=-\f d{24},\ b_3'=-\f{d(d+1)}2\in\Z,
\\b_4'=&108 d + 343 d^2  + 108 d^3,
\\ b_5'=&720000 d + 4663750 d^2  + 4663750 d^3  + 720000 d^4.
\end{align*}
\end{remark}

Let $n\in\Z^+=\{1,2,3,\ldots\}$. For any polynomial $P(z)=\sum_{k=0}^{n-1}a_kz^k$ of degree $n-1$ with complex coefficients, it is known (cf. \cite[Lemma 9]{K05}) that
\begin{equation}\label{P}|P(x_i+y_j)|_{1\le  i,j\le n}=a_{n-1}^n\prod_{k=0}^{n-1}\bi{n-1}k\times\prod_{1\le i<j\le n}(x_i-x_j)(y_j-y_i).
\end{equation}
(Note that \cite[Lemma 9]{K05} mis-wrote $\prod_{k=0}^{n-1}\bi{n-1}k$ as $\prod_{k=1}^n\bi{n-1}k$ which is zero.) By letting $a_{n-1}\to 0$ we see that \eqref{P} also holds if $\deg P<n-1$.
In particular,
$$\l|(i+j)^k\r|_{0\le i,j\le n-1}=0\quad\t{for all}\ k=0,\ldots,n-2,$$
and
\begin{align*}\l|(i+j)^{n-1}\r|_{0\le i,j\le n-1}=&(-1)^{n(n-1)/2}\prod_{0\le i<j\le n-1}(j-i)^2\times\prod_{k=0}^{n-1}\bi {n-1}k
\\=&(-1)^{n(n-1)/2}((n-1)!)^n.
\end{align*}
But this is of no help in simplifying the determinants $a_n$ and $b_n$ given in \eqref{ab}
even if $c=0$ and $d=1$.

Our third theorem confirms Conjecture 4.3 of Sun \cite{S19}.

\begin{theorem}\label{Th1.3} Let $p$ be an odd prime, and let $c,d\in\Z$ with $p\nmid cd$. Define
$$S_c(d,p):=\bg|\l(\f{i^2+dj^2+c}p\r)\bg|_{1\le i,j\le(p-1)/2}.$$
Then
\begin{equation}\label{Sc}\l(\f {S_c(d,p)}p\r)=\begin{cases}1&\t{if}\ (\f cp)=1\ \t{and}\ (\f dp)=-1,
\\(\f{-1}p)&\t{if}\ (\f cp)=(\f dp)=-1,
\\(\f{-2}p)&\t{if}\ (\f{-c}p)=(\f dp)=1,
\\(\f{-6}p)&\t{if}\ (\f{-c}p)=-1\ \t{and}\ (\f dp)=1.\end{cases}
\end{equation}
\end{theorem}

In contrast, for any odd prime $p$ and $d\in\Z$ with $p\nmid d$, Sun \cite[(1.15) and (1.20)]{S19} showed that
$$\l(\f{S_0(d,p)}p\r)=\begin{cases}(\f{-1}p)&\t{if}\ (\f dp)=1,
\\0&\t{if}\ (\f dp)=-1.\end{cases}$$
But the method used to prove this does not work for Theorem \ref{Th1.3}.

We will prove Theorem \ref{Th1.1} in the next section.
Using an auxiliary formula in Section \ref{sect.aux}, we are going to prove Theorems \ref{Th1.2}
and \ref{Th1.3}
in Sections \ref{sect.pf2} and \ref{Sect.5} respectively.

\section{\label{sect.pf1}Proof of Theorem \ref{Th1.1}}
\setcounter{lemma}{0}
\setcounter{theorem}{0}
\setcounter{corollary}{0}
\setcounter{remark}{0}
\setcounter{equation}{0}

\begin{lemma}\label{Lem2.1} Let $p$ be a prime and let $k\in\N=\{0,1,2,\ldots\}$. Then
$$\sum_{i=1}^{p-1}i^k\eq\begin{cases}-1\pmod p&\t{if}\ p-1\mid k,
\\0\pmod p&\t{if}\ p-1\nmid k.\end{cases}$$
\end{lemma}

This is a well known fact, see, e.g., \cite[Section 15.2, Lemma 2]{IR}.

For a prime $p$, any rational number with its denominator not divisible by $p$
is a $p$-adic integer. We deal with $p$-adic congruences in our proof of Theorem \ref{Th1.1} and further on.

\medskip
\noindent{\it Proof of Theorem \ref{Th1.1}}.
If $n$ is composite, then $n$ can be written as $n = pm$ for some odd integers $m \ge p > 2$, and
thus $i := (m-p)/2$ and $i' := (m+p)/2$ are integers satisfying
$0 \leq i < i' \le (n-1)/2$ and $i^2 \equiv i'^2\pmod n$.
So, when $n$ is composite, there are $0\le i<i'\le (n-1)/2$ such that
$$(i^2+dj^2)\l(\f{i^2+dj^2}n\r)\eq ((i')^2+dj^2)\l(\f{(i')^2+dj^2}n\r)\pmod n$$
for all $j=0,\ldots,(n-1)/2$,
and hence \eqref{0(mod n)} holds
(since an integer matrix that has two rows congruent to each other modulo $n$
must have a determinant congruent to $0$ modulo $n$).

It thus remains to prove Theorem \ref{Th1.1} in the case when
$n$ is a prime.
So let us assume that $n$ is a prime $p>3$.

Fix $j \in \left\{0,\ldots,(p-1)/2\right\}$. As
$(\f ap) \equiv a^{(p-1)/2} \pmod p$ for all $a\in\Z$ (due to Euler), we have
\begin{align*}
&\sum_{i=1}^{(p-1)/2}(i^2+dj^2)\l(\f{i^2+dj^2}p\r)
\\
\eq&\sum_{i=1}^{(p-1)/2}(i^2+dj^2)^{(p+1)/2}
=\sum_{i=1}^{(p-1)/2}\ \ \sum_{k=0}^{(p+1)/2}\bi{(p+1)/2}ki^{2k}(dj^2)^{(p+1)/2-k}
\\
=&\sum_{k=0}^{(p+1)/2}\bi{(p+1)/2}k(dj^2)^{(p+1)/2-k}\sum_{i=1}^{(p-1)/2}i^{2k}\pmod p.
\end{align*}
Multiplying this by $2$, we obtain
\begin{align}
&2 \sum_{i=1}^{(p-1)/2}(i^2+dj^2)\l(\f{i^2+dj^2}p\r)\nonumber
\\
\eq&\sum_{k=0}^{(p+1)/2}\bi{(p+1)/2}k(dj^2)^{(p+1)/2-k}\sum_{i=1}^{(p-1)/2}(i^{2k}+(p-i)^{2k})
\nonumber\\
\eq&\sum_{k=0}^{(p+1)/2}\bi{(p+1)/2}k(dj^2)^{(p+1)/2-k}\sum_{i=1}^{p-1}i^{2k}\pmod p.
\label{pf.Th1.1.3}
\end{align}
For each $k\in\{0,\ldots,(p+1)/2\}$, clearly
$$p-1\mid 2k\iff k=0\ \t{or}\ k=(p-1)/2$$
as $p>3$, and hence by Lemma \ref{Lem2.1} we get
\[
\sum_{i=1}^{p-1} i^{2k} \eq\begin{cases}-1\pmod p&\t{if}\ k=0\ \t{or}\ k=(p-1)/2,
\\0\pmod p&\t{otherwise}.\end{cases}
\]
Thus, \eqref{pf.Th1.1.3} simplifies to
\begin{align*}&2\sum_{i=1}^{(p-1)/2}(i^2+dj^2)\l(\f{i^2+dj^2}p\r)
\\\eq&\bi{(p+1)/2}0(dj^2)^{(p+1)/2}(-1)+\bi{(p+1)/2}{(p-1)/2}(dj^2)(-1)
\\\eq&-dj^2\l(\l(\f{dj^2}p\r)+\f{p+1}2\r)\eq-\f{dj^2}2\l(\f {dj^2}p\r)\l(2+\l(\f dp\r)\r)
\pmod p
\end{align*}
(where the last two congruence signs relied on
$(\f ap) \equiv a^{(p-1)/2} \pmod p$ and on
the easily verified congruence
$dj^2 \equiv dj^2 (\f {dj^2}p) (\f dp) \pmod p$,
respectively).
Note that $2+(\f dp)$ is relatively prime to $p$ since $p>3$. Therefore
$$\sum_{i=1}^{(p-1)/2}\f 4{2+(\f dp)}(i^2+dj^2)\l(\f{i^2+dj^2}p\r)+(0^2+dj^2)\l(\f{0^2+dj^2}p\r)\eq0\pmod p .$$

The last congruence holds for all $j=0,\ldots,(p-1)/2$.
Thus, if we add the last $(p-1)/2$ rows multiplied by $4/(2+(\f dp))$ to the first row
of the determinant
$$D_p:=\l|(i^2+dj^2)\l(\f{i^2+dj^2}p\r)\r|_{0\le i,j\le (p-1)/2},$$
then all the entries in the first row of the resulting determinant are multiples of $p$.
So we have $D_p\eq0\pmod p$ as desired.

In view of the above, this completes the proof of Theorem \ref{Th1.1}. \qed

\section{\label{sect.aux}A general formula for $\l|(x_i+y_j)^n\r|_{1\le i,j\le n}$}
\setcounter{lemma}{0}
\setcounter{theorem}{0}
\setcounter{corollary}{0}
\setcounter{remark}{0}
\setcounter{equation}{0}

For each $k=1,\ldots,n$, the $k$th elementary symmetric polynomial $\sigma_k$
in $x_1,\ldots,x_n$ is defined by
$$\sigma_k(x_1,\ldots,x_n)=\sum_{1\le i_1<\ldots<i_k\le n}\prod_{j=1}^kx_{i_j}.$$
In addition, we set $\sigma_0(x_1,\ldots,x_n)=1$ as usual.

To prove Theorem \ref{Th1.2}, we need the following auxiliary theorem
which improves a result of \cite[\S354(a), pp. 349--350]{MM}.

\begin{theorem}\label{Th3.1} Let $n$ be a positive integer, and let $x_1, \ldots, x_n, y_1, \ldots, y_n$
be elements of any commutative ring with identity. Then
\begin{equation}\label{general}\begin{aligned}
\l|(x_i+y_j)^n\r|_{1\le i,j\le n}
= \ &(-1)^{n(n-1)/2}\prod_{1\le i<j\le n}(x_j-x_i)(y_j-y_i)
\\
\times\sum_{k=0}^n&\(\prod_{r \in [0,n] \sm \{k\}} \bi nr\)
\sigma_k(x_1,\ldots,x_n)\sigma_{n-k}(y_1,\ldots,y_n),
\end{aligned}\end{equation}
where $[0,n]$ denotes the set $\{0,\ldots,n\}$.
\end{theorem}

\Proof. Define an $n\times\left(n+1\right)$-matrix $A$ and an
$\left(n+1\right)\times n$-matrix $B$ by
$$A=\l[\bi nkx_i^k\r]_{1\le i\le n\atop 0\le k\le n}
\ \ \t{and}\ \ B=\l[y_j^{n-k}\r]_{0\le k\le n\atop 1\le j\le n}.$$
As the binomial formula yields
$$(x_i+y_j)^n=\sum_{k=0}^n\bi nk x_i^ky_j^{n-k},$$
we have
$$AB=\l[(x_i+y_j)^n\r]_{1\le i,j\le n}.$$
Applying the Cauchy-Binet formula, we therefore get
\begin{align}
&\l|(x_i+y_j)^n\r|_{1\le i,j\le n}
\\=&\sum_{k=0}^n\l|\bi njx_i^j\r|_{1\le i\le n\atop j\in[0,n]\sm\{k\}}\l|y_j^{n-i}\r|_{1\le j\le n\atop i\in[0,n]\sm\{k\}}
\nonumber\\
=&\sum_{k=0}^n \(\prod_{r \in [0,n] \sm \{k\}} \bi nr\)
\l|x_j^i\r|_{i\in[0,n]\sm\{k\}\atop 1\le j\le n}\times
(-1)^{n(n-1)/2}\l|y_j^i\r|_{i\in[0,n]\sm\{n-k\}\atop 1\le j\le n}
\nonumber\\&
\label{pf.Th2.1.4}
\end{align}
(by standard properties of determinants).
For each $k\in\{0,\ldots,n\}$, comparing the coefficient of $x^k$ on
both sides of the polynomial equality
$$\left|\begin{matrix}
1&1&\cdots &1&1\\x_1&x_2&\cdots&x_{n}&x\\
\vdots&\vdots&\ddots&\vdots&\vdots\\
x_1^{n-1}&x_2^{n-1}&\cdots&x_{n}^{n-1}&x^{n-1}
\\x_1^n&x_2^n&\cdots&x_n^n&x^n
\end{matrix}\right|=\prod_{1\le i<j\le n}(x_j-x_i)\times\prod_{i=1}^n(x-x_i)$$
(a consequence of Vandermonde's determinant),
we find that
$$(-1)^{n-k}\l|x_j^i\r|_{i\in[0,n]\sm\{k\}\atop 1\le j\le n}=(-1)^{n-k}\sigma_{n-k}(x_1,\ldots,x_n)
\prod_{1\le i<j\le n}(x_j-x_i)$$
(where the left-hand side was computed by expanding
the determinant along its last column).
Hence,
$$\l|x_j^i\r|_{i\in[0,n]\sm\{k\}\atop 1\le j\le n}=\sigma_{n-k}(x_1,\ldots,x_n)
\prod_{1\le i<j\le n}(x_j-x_i).$$
Similarly,
$$\l|y_j^i\r|_{i\in[0,n]\sm\{n-k\}\atop 1\le j\le n}=\sigma_{k}(y_1,\ldots,y_n)
\prod_{1\le i<j\le n}(y_j-y_i).$$
Therefore, we can rewrite \eqref{pf.Th2.1.4} as
\begin{align*}\l|(x_i+y_j)^n\r|_{1\le i,j\le n}=&(-1)^{n(n-1)/2}\prod_{1\le i<j\le n}(x_j-x_i)(y_j-y_i)
\\\times\sum_{k=0}^n&\(\prod_{r \in [0,n] \sm \{k\}} \bi nr \)\sigma_{n-k}(x_1,\ldots,x_n)\sigma_k(y_1,\ldots,y_n).
\end{align*}
Substituting $n-k$ for $k$ on the right-hand side, and observing that
$\prod_{r \in [0,n] \sm \{n-k\}} \bi nr = \prod_{r \in [0,n] \sm \{k\}} \bi nr$,
we obtain the desired \eqref{general}. \qed

\section{\label{sect.pf2}Proof of Theorem \ref{Th1.2}}
\setcounter{lemma}{0}
\setcounter{theorem}{0}
\setcounter{corollary}{0}
\setcounter{remark}{0}
\setcounter{equation}{0}

\medskip
\noindent{\it Proof of Theorem \ref{Th1.2}}. Clearly \eqref{divide} holds if $a_n',b_n'\in\Z$.

(i) Let us first discuss $a_n$ and $a_n'\in\Z$.
By Remark \ref{Rem1.2} we have the desired result for $a_n$ and $a_n'$
when $n=3,4$; so let us assume that $n \geq 5$.

Define
$$S(n):= \sum_{k=0}^n\sigma_k(0,\ldots,n-1)\sigma_{n-k}(d0+c,\ldots,d(n-1)+c) \prod_{r \in [0,n] \sm \{k\}} \bi nr,$$
which is positive if $c\ge0$ and $d>0$.
Applying Theorem \ref{Th3.1}, we find that
\begin{align*}a_n=&(-1)^{n(n-1)/2}\prod_{0\le i<j\le n-1}(j-i)(dj-di)\times S(n)
\\=&(-d)^{n(n-1)/2}S(n)\prod_{j=1}^{n-1}(j!)^2.
\end{align*}
Hence, $(-1)^{n(n-1)/2}a_n>0$ if $c\ge0$ and $d>0$.
To prove $a_n'\in\Z$ it suffices to show that
\begin{equation}\label{n^2}1!2!\ldots(n-3)!S(n)\eq0\pmod{n^2}.
\end{equation}

Fix $k \in [0, n]$.
The product $\prod_{r \in [0,n] \sm \{k\}} \bi nr$ contains
at least two of the three factors $\bi n1$, $\bi n2$ and
$\bi n{n-1}$ (since $n \geq 5$).
But each of these three factors is divisible by $n$ or
(in the case of the second factor) by $n/2$ (when $n$ is even).
Thus, the product $\prod_{r \in [0,n] \sm \{k\}} \bi nr$ is
divisible by $n \times n$ or (when $n$ is even) by $n\times n/2$.
In either case, it follows that
$2\prod_{r \in [0,n] \sm \{k\}} \bi nr$ is divisible by
$n^2$.
Since $n \geq 5$ yields $2 \mid 1!2!\ldots(n-3)!$, we thus conclude
that $1!2!\ldots(n-3)! \prod_{r \in [0,n] \sm \{k\}} \bi nr$
is divisible by $n^2$.
Since we have shown this for all $k \in [0, n]$, it follows
that $1!2!\ldots(n-3)!S(n)$ is divisible by $n^2$.
This proves \eqref{n^2}.
\medskip

(ii) For $n=3,4,5$, we have the desired result for $b_n$ and $b_n'$ by Remark \ref{Rem1.2}.

Now we assume $n\ge6$ and define
\begin{align*}T(n):=&\sum_{k=0}^n\sigma_k(0^2,\ldots,(n-1)^2)\sigma_{n-k}(d0^2,\ldots,d(n-1)^2) \prod_{r \in [0,n] \sm \{k\}} \bi nr
\\=&\sum_{k=1}^{n-1}d^{n-k}\sigma_k(0^2,\ldots,(n-1)^2)\sigma_{n-k}(0^2,\ldots,(n-1)^2) \prod_{r \in [0,n] \sm \{k\}} \bi nr,
\end{align*}
which is positive if $d>0$
and always satisfies $d \mid S(n)$.
In view of Theorem \ref{Th3.1},
\begin{align*}
\l|(i^2+dj^2)^n\r|_{0\le i,j\le n-1}
=& (-1)^{n(n-1)/2}\prod_{0\le i<j\le n-1}(j^2-i^2)(dj^2-di^2)\times T(n) \\
=& (-d)^{n(n-1)/2}T(n)\prod_{j=1}^{n-1}((2j-1)!j)^2
\end{align*}
and so
\begin{equation}\label{bn} b_n=(-d)^{n(n-1)/2}T(n)((n-1)!)^2\prod_{j=1}^{n-1}((2j-1)!)^2.
\end{equation}
Thus $(-1)^{n(n-1)/2}b_n>0$ if $d>0$.
Also, \eqref{bn} yields
\begin{align}
(-1)^{n(n-1)/2}b_n'
&=\f{((n-1)!)^2\prod_{j=1}^{n-1}(2j-1)!}{(2n-1)!\times 2\prod_{k=1}^nk!}T(n)
\nonumber\\
&=\f{\prod_{k=1}^{n-2}\f{(2k-1)!}{k!}}{2n(2n-1)(2n-2)}T(n).
\label{T(n)}
\end{align}

Note that
$$\f1{4n(n-1)}\prod_{k=1}^{n-2}\f{(2k-1)!}{k!}=\begin{cases}105&\t{if}\ n=6,\\226800&\t{if}\ n=7,
\\9430344000&\t{if}\ n=8.\end{cases}$$
If $n\ge 9$, then  $n+4\le 2n-5$ and one of $n+1,n+2,n+3,n+4$ is divisible by $4$, hence
$$4n(n-1)\ \bg|\ \f{(2(n-2)-1)!}{(n-2)!}.$$
So we always have
$$2n(2n-2)\ \bigg|\ \prod_{k=1}^{n-2}\f{(2k-1)!}{k!}.$$
Therefore, \eqref{T(n)} leads to
$(2n-1) b_n' \in \Z$.
If we can furthermore show that
$2n (2n-2) b_n' \in \Z$, then we will conclude that
$b_n' \in \Z$ (since $2n-1$ is coprime to
$2n (2n-2)$).
So it remains to show that $2n (2n-2) b_n' \in \Z$, i.e., that
$2n-1 \mid T(n)\prod_{k=1}^{n-2}\f{(2k-1)!}{k!} $.

If $2n-1=pq$ with $p,q\in\Z^+$ and $3\le p<q$, then
$p<q\le \f{2n-1}3\le n-3$, and hence $2n-1=pq$ divides  $$\f{(2(p-1)-1)!}{(p-1)!}\times\f{(2(q-1)-1)!}{(q-1)!}.$$
If $2n-1=p^2$ with $p$ an odd prime, then $5\le p=\sqrt{2n-1}<n-2$ since $2n-1>n\ge9$, hence $2n-1=p^2$ divides
$$\f{(2(p-2)-1)!}{(p-2)!}\times \f{(2(p-1)-1)!}{(p-1)!}.$$
If $2n-1$ is a prime $p$, then $p>3$ and
\begin{align*}
b_n
&= \l|(i^2+dj^2)^{(p+1)/2}\r|_{0\le i,j\le (p-1)/2}
\\
&\eq\l|(i^2+dj^2)\l(\f{i^2+dj^2}p\r)\r|_{0\le i,j\le (p-1)/2}
\\
&\eq0\pmod p
\end{align*}
by Theorem \ref{Th1.1}, hence $2n-1=p$ divides $T(n)$ by \eqref{bn} and due to $d\mid T(n)$.
In either case, we obtain $2n-1 \mid T(n)\prod_{k=1}^{n-2}\f{(2k-1)!}{k!}$.

\medskip
The proof of Theorem \ref{Th1.2} is now complete. \qed
\medskip

\section{\label{Sect.5}Proof of Theorem \ref{Th1.3}}
\setcounter{lemma}{0}
\setcounter{theorem}{0}
\setcounter{corollary}{0}
\setcounter{remark}{0}
\setcounter{equation}{0}

We need the following known lemma (see \cite[Lemma 2.3]{S19}):

\begin{lemma} \label{Lem5.2new}
Let $p$ be a prime with $p \equiv 1 \pmod 4$, and write $n = (p-1)/2$.
Then
\[
\l( \f {n!}p \r) = \l( \f 2p \r) .
\]

\end{lemma}

\medskip
\noindent{\it Proof of Theorem \ref{Th1.3}}. For convenience we set $n=(p-1)/2$.
Applying Theorem \ref{Th3.1}, we see that
\begin{align}
& \l|(i^2+dj^2+c)^{n}\r|_{1\le i,j\le n}
\nonumber\\
=&(-1)^{n(n-1)/2}\prod_{1\le i<j\le n}(j^2-i^2)(dj^2+c-(di^2+c))\times\prod_{r=0}^n\bi nr
\nonumber\\
&\times\sum_{k=0}^n\f{\sigma_k(1^2,\ldots,n^2)\sigma_{n-k}(d1^2+c,\ldots,dn^2+c)}{\bi nk}
\nonumber\\
=&(-d)^{n(n-1)/2}\prod_{1\le i<j\le n}(j^2-i^2)^2\times\prod_{r=0}^n\bi nr\times R_n,
\label{pf.Th1.3.0}
\end{align}
where
\begin{align}
R_n:=&\sigma_n(1^2,\ldots,n^2)+\sigma_n(d1^2+c,\ldots,dn^2+c) \nonumber\\
&+\sum_{0<k<n}\f{\sigma_k(1^2,\ldots,n^2)\sigma_{n-k}(d1^2+c,\ldots,dn^2+c)}{\bi nk}.
\label{pf.Th1.3.1}
\end{align}

As observed in \cite[(3.2)]{S19}, we have the polynomial congruence
\begin{align}
 & \sum_{k=0}^n(-1)^k\sigma_k(1^2,\ldots,n^2)x^{n-k} \nonumber\\
= &\prod_{r=1}^n(x-r^2)\eq x^n-1\pmod p.
\label{pf.Th1.3.2}
\end{align}
So $\sigma_n(1^2,\ldots,n^2)\eq -(-1)^n\pmod p$ and $\sigma_k(1^2,\ldots,n^2)\eq0\pmod p$
for all $k=1,\ldots,n-1$. Note also that $p\nmid \bi nk$ for all $k=0,\ldots,n$. Therefore,
\eqref{pf.Th1.3.1} yields
\begin{align*}&R_n+(-1)^n
\\\eq&\sigma_n(d1^2+c,\ldots,dn^2+c)=\prod_{r=1}^n(c+dr^2)
=(-d)^n\prod_{r=1}^n\l(-\f cd-r^2\r)
\\\eq&(-d)^n\l(\l(-\f cd\r)^n-1\r)=c^n-(-d)^n\eq\l(\f cp\r)-\l(\f{-d}p\r)
\pmod p ,
\end{align*}
where we have used \eqref{pf.Th1.3.2} in the third-to-last
congruence.
Solving this for $R_n$ and substituting the result
into \eqref{pf.Th1.3.0}, and noting that the left-hand side of \eqref{pf.Th1.3.0}
is congruent to $S_c(d, p)$ modulo $p$ (since $( \f ap ) \equiv a^n \pmod p$ for all $a \in \Z$),
we find
\begin{align}
&S_c(d,p) \nonumber\\
\eq&(-d)^{n(n-1)/2}\l(\prod_{1\le i<j\le n}(j^2-i^2)^2\r) \times \prod_{r=0}^n\bi nr
\nonumber\\
&\quad \times \l(\l(\f cp\r)-\l(\f{-d}p\r)-(-1)^n\r)
\pmod p.
\label{pf.Th1.3.5}
\end{align}

Clearly,
\[
\prod_{r=0}^n\bi nr = \dfrac{(n!)^{n+1}}{\l( 0! 1! \cdots n! \r)^2}
\]
and hence
\begin{align*}
\l(\f{\prod_{r=0}^n\bi nr}p\r)
= \l(\f {(n!)^{n+1}}p\r)
= \l(\f {n!}p\r)^{n+1}
= \l(\f 2p\r)^{n+1}
\end{align*}
(by using Lemma~\ref{Lem5.2new} when $2\mid n$).
In view of the known identity
$( \f 2p ) =(-1)^{(p^2-1)/8}= (-1)^{n(n+1)/2},$
we can rewrite this as
\[
\l(\f{\prod_{r=0}^n\bi nr}p\r)
= ((-1)^{n(n+1)/2})^{n+1}
= (-1)^{n(n+1)^2/2} .
\]
Hence, \eqref{pf.Th1.3.5} yields
\begin{equation}\label{Legen}
\l(\f{S_c(d,p)}p\r)=\l(\f{-d}p\r)^{n(n-1)/2} (-1)^{n(n+1)^2/2} \l(\f{(\f cp)-(\f{-d}p)-(-1)^n}p\r).
\end{equation}
Note that $(\f{-1}p)=(-1)^n$ and
$$(-1)^{n(n+1)^2/2-n^2(n-1)/2}=(-1)^{n(3n+1)/2}=(-1)^{n(n-1)/2}.$$
So \eqref{Legen} can be rewritten as
\begin{equation}\label{Legen'}
\l(\f{S_c(d,p)}p\r)= (-1)^{n(n-1)/2} \l(\f{d}p\r)^{n(n-1)/2}\l(\f{(\f cp)-(-1)^n(1+(\f dp))}p\r).
\end{equation}

Now it remains to deduce \eqref{Sc} from \eqref{Legen'}.
\medskip

{\it Case} 1. $(\f cp)=1$ and $(\f dp)=-1$.

In this case,
\eqref{Legen'} becomes
\begin{align*}
\l(\f{S_c(d,p)}p\r)=(-1)^{n(n-1)/2}(-1)^{n(n-1)/2} \l(\f{1-(-1)^n(1-1)}p\r)=1.
\end{align*}

{\it Case} 2. $(\f cp)=(\f dp)=-1$.

In this case, \eqref{Legen'} gives
\begin{align*}
\l(\f{S_c(d,p)}p\r)
&=(-1)^{n(n-1)/2} (-1)^{n(n-1)/2} \l(\f {-1-(-1)^n(1-1)}p\r)
= \l(\f {-1}p\r).
\end{align*}

{\it Case} 3. $(\f {-c}p)=(\f dp)=1$.

In this case,  $(\f cp) =(\f{-1}p)=(-1)^n$.
Hence, \eqref{Legen'} yields
\begin{align*}
\l(\f{S_c(d,p)}p\r)
&= (-1)^{n(n-1)/2}\l(\f{(-1)^n - (-1)^n2}p\r)=(-1)^{n(n+1)/2-n}\l(\f{-1}p\r)^{n+1}
\\
&= (-1)^{(p^2-1)/8}(-1)^n(-1)^{n(n+1)}=\l(\f 2p\r)\l(\f{-1}p\r)
= \l(\f {-2}{p} \r) .
\end{align*}

{\it Case} 4. $(\f {-c}p)=-1$ and $(\f dp)=1$.

In this case, $(\f cp) =-(\f{-1}p)=-(-1)^n$. Hence, by \eqref{Legen'} we have
\begin{align*}
\l(\f{S_c(d,p)}p\r)
&= (-1)^{n(n-1)/2} \l(\f{-(-1)^n - (-1)^n2}p\r)
\\
&= (-1)^{n(n+1)/2-n} \l(\f{-1}p\r)^n\l(\f{-3}p\r)
\\&= (-1)^{n(n+1)/2} (-1)^n \l((-1)^n\r)^n\l(\f{-3}p\r) \\
&= (-1)^{n(n+1)/2} \l(\f{-3}p\r)
=\l(\f2p\r)\l(\f{-3}p\r)= \l(\f {-6}{p} \r).
\end{align*}

In view of the above, \eqref{Sc} holds as desired. This concludes the proof. \qed

\Ack. We thank Prof. Guo-Niu Han and the anonymous referee for helpful comments.

\setcounter{conjecture}{0} \end{document}